\documentclass[a4,psfig,english,11pt,epsf,portrait]{article}
\usepackage{amsmath,amsfonts,latexsym,amscd,amssymb,theorem}
\usepackage{epsfig}
\usepackage{amsmath}
\usepackage{color}
\definecolor{purple}{rgb}{0.65, 0, 1}
\definecolor{orange}{rgb}{1,.5,0}
\definecolor{brown}{rgb}{.9,.73,.26}

\def\<{\langle}
\def\>{\rangle}
\newcommand{\ba}{\begin{eqnarray}}
\newcommand{\ea}{\end{eqnarray}}

\newtheorem{thm}{Theorem}[section]

\newtheorem{Theorem}[thm]{Theorem}
\newtheorem{Definition}[thm]{Definition}
\newtheorem{Lemma}[thm]{Lemma}
\newtheorem{Proposition}[thm]{Proposition}

\newtheorem{Remark}[thm]{Remark}

\numberwithin{equation}{section}

\begin{document}
\title{\bf Nonexistence of global positive solutions for p-Laplacian equations with non-linear memory}
\author{
\normalsize\textsc{ Mokhtar Kirane\footnote{Department of Mathematics,  Faculty of Arts and Science,  Khalifa University, P.O. Box: 127788,  Abu Dhabi, UAE; mokhtar.kirane@ku.ac.ae.
\newline \indent $\,\,{}^{a}$Department of Mathematics, Sultan Qaboos University,  FracDiff Research Group (DR/RG/03),  P.O. Box 46, Al-Khoud 123, Muscat, Oman; ahmad.fino01@gmail.com; a.fino@squ.edu.om.
\newline \indent $\,\,{}^{b}$Department of Mathematics, Sultan Qaboos University,  FracDiff Research Group (DR/RG/03),  P.O. Box 36, Al-Khoud 123, Muscat, Oman; skerbal@squ.edu.om.}, Ahmad Z. Fino$^{a}$, Sebti Kerbal$^{b}$
  }}

\vspace{20pt}

\maketitle

%%%%%%%%%%%%%%%%%%%%%%%%%%%%%%%%%%%%%%%%%%%%%%%%%%%%%%%%%%%%%%%%%%%%%%%%%%
%%%%%%%%%%%%%%%%%%%%%%%%%%%%%%%%%%%%%%%%%%%%%%%%%%%%%%%%%%%%%%%%%%%%%%%%%%

\centerline{\small{\bf{Abstract}}} \noindent{\small{The Cauchy problem in $\mathbb{R}^d,$ $d\geq 1,$ for a non-local in time p-Laplacian equations is considered. The nonexistence of nontrivial global weak solutions by using the test function method is obtained.}} \hfill\break


\section{Introduction}

In this paper, we consider the following problem
\begin{equation}\label{eq1}
\left\{\begin{array}{ll}
\,\, \displaystyle {u_{t}-\mbox{div}\left(|\nabla u|^{p-2}\nabla u\right) =\mathcal{I},} &\displaystyle {t>0,x\in {\mathbb{R}^d},}\\
\displaystyle u\geq0, &\displaystyle {t>0,x\in {\mathbb{R}^d},}\\
\displaystyle{u(0,x)= u_0(x),\qquad\qquad}&\displaystyle{x\in {\mathbb{R}^d},}
\end{array}
\right.
\end{equation}
where $\mathcal{I}:=  \int_0^t(t-s)^{-\gamma}u(s)^{q}\,ds,$ $p,q>1$, $d\geq1$, $0<\gamma<1$, and $u_0\in L_{loc}^2(\mathbb{R}^d)$.  We are interested in  the nonexistence of nontrivial global weak solutions.

The study of nonexistence of global solutions for nonlinear parabolic equations started by Fujita \cite{Fujita}; he studied the  Cauchy
problem
\begin{equation}\label{eq2}
\left\{\begin{array}{ll}
\,\, \displaystyle {u_{t}-\Delta u =u^p,} &\displaystyle {p >1, t>0,x\in {\mathbb{R}^d},}\\
{}\\
\displaystyle{u(0,x)= u_0(x)\geq0,\qquad\qquad}&\displaystyle{x\in {\mathbb{R}^d};}
\end{array}
\right.
\end{equation}
and obtained that

$a)$  there is no nontrivial global solution of (\ref{eq2}) whenever $p< 1 + 2/d$. 

$b)$  whereas  (\ref{eq2}) admits a global solution whenever $p > 1 + 2/d$ and $u_0(x)\leq\delta e^{-|x|^2}$, ($0 <\delta\ll 1$),.\\
For the limiting case $p = 1 + 2/d$,  it is shown by Hayakawa \cite{Hayakawa} for $d = 1, 2$ and by Kobayashi, Sirao and Tanaka \cite{kobayashi} for $d \geq 1$ that (\ref{eq2})
has no global solution $u(x, t)$ satisfying $\|u(\cdotp, t)\|_{\infty} < \infty$ for $t\geq 0$.  Weissler in \cite{weissler} proved that if (\ref{eq2})
has no global solution $u(x, t)$ satisfying $\|u(\cdotp, t)\|_{q} < \infty$ for $t> 0$ and some $q\in [1,+\infty)$ whenever  if $p = 1+2/d$.  Therefore,  the value
\begin{equation}\label{p_F}
p_F = 1 + 2/d
\end{equation}
 is the limiting exponent of (\ref{eq2}); it is called Fujita's exponent.  

Since Fujita's paper,  a sizeable number of extensions in many directions have been published.  Recently,  Cazenave,  Dickstein and Weissler \cite{CDW} extended the results of Fujita 
to the non-local in time heat  equation 
\begin{equation}\label{eq3}
\left\{\begin{array}{ll}
\,\, \displaystyle {w_t-\Delta w =\int_0^t(t-s)^{-\vartheta}|w(s)|^{p-1}w(s)\,ds,} &\displaystyle {t>0,x\in {\mathbb{R}^d}, \; p>1,}\\
{}\\
\displaystyle{w(0,x)= w_0(x),\qquad\qquad}&\displaystyle{x\in {\mathbb{R}^d},}
\end{array}
\right.
\end{equation}
where  $w_0\in C_0(\mathbb{R}^d)$ and $\vartheta\in(0,1)$. In this case, the value of the critical Fujita exponent is
\begin{equation}\label{p_c}
p_c=\max\left\{\frac{1}{\vartheta},1+\frac{2(2-\vartheta)}{(d-2(1-\vartheta))_+}\right\}.
\end{equation}

For the p-Laplacian equation
\begin{equation}\label{eq4}
\left\{\begin{array}{ll}
\displaystyle {w_{t}-\mbox{div}\left(|\nabla w|^{p-2}\nabla w\right) =w^{q},} &\displaystyle {t>0,x\in {\mathbb{R}^d}, p>2, q>1,}\\
\displaystyle w\geq0, &\displaystyle {t>0,x\in {\mathbb{R}^d},}\\
\displaystyle{w(0,x)= w_0(x)\geq0,\qquad\qquad}&\displaystyle{x\in {\mathbb{R}^d},}
\end{array}
\right.
\end{equation}
Zhao \cite{zhao} and  Mitidieri and Pohozaev \cite{PM1} obtained the critical exponent
\begin{equation}\label{q^*}
q^* = p -1 + p/d;
\end{equation}
in fact,  Zhao \cite{zhao} proved that if $p-1< q < q^*$,  then the Cauchy problem (\ref{eq4}) has no nontrivial global solution,  however if $q > q^*$ and $w_0(x)$ is small enough,  then (\ref{eq4}) admits a global solution.  Mitidieri and Pohozaev \cite{PM1} completed the study by proving the nonexistence of nontrivial global solution in the case $q\leq q^*$ and for all $p>2d/(d+1)$.  Andreucci and Tedeev \cite{andreucci, liu} obtained similar results by considering doubly singular parabolic equations.

The test function method was used to prove the nonexistence of global solutions. This method was introduced by  Baras and Kersner in \cite{Baras} and developed by Zhang in \cite{Zhang} and Pohozaev and Mitidieri in \cite{PM1},  it was also used by Kirane et al.  in \cite{KLT}.  

Here,  we are concerned with the non-existence of nontrivial global solutions of (\ref{eq1});   inspired by \cite{PM1},  we choose a suitable test function in the weak formulation of the problem.  Our main result is
\begin{Theorem}\label{theo2}${}$
Let $u_0\in L^2_{loc}(\mathbb{R}^d),$ $0<\gamma<1$, $p>1$, $d\geq1$,  $q>\max\{1,p-1\}$. Suppose that
\begin{equation}\label{C0}
\left\{\begin{array}{ll}
\displaystyle {q\leq\max\left\{\frac{p-1}{\gamma};q_c\right\},} &\displaystyle {\hbox{if}\,\, p<\frac{d}{1-\gamma},}\\
\displaystyle {q<\infty,} &\displaystyle {\hbox{if}\,\, p\geq\frac{d}{1-\gamma},}\\
\end{array}
\right.
\end{equation}
where
$$q_c:=p-1+\frac{(d+p)(1-\gamma)(p-1)+p-d(1-\gamma)}{d-(1-\gamma)p},$$
then  no nontrivial global weak solutions exist for problem $(\ref{eq1})$.
\end{Theorem}
\begin{Remark}
Note that (\ref{C0}) can be seen as follows
$$
q\leq\max\left\{\frac{p-1}{\gamma};p-1+\frac{(d+p)(1-\gamma)(p-1)+p-d(1-\gamma)}{(d-(1-\gamma)p)_+}\right\},$$
whence\\
$\bullet$ when $p=2$ and $\gamma\rightarrow 1$, then $q_c\rightarrow p_F$ where $p_F$ is defined in (\ref{p_F}).\\
$\bullet$ when $p=2$, then $q_c=p_c$ where $p_c$ is defined in (\ref{p_c}).\\
$\bullet$ when $\gamma\rightarrow 1$, and $p>2d/(d+1)$, then $q_c\rightarrow q^*$ where $q^*$ is defined in (\ref{q^*}).\\
\end{Remark}
\begin{Remark}
The same result can be obtained for the problem
\[
\left\{\begin{array}{ll}
\displaystyle {u_{t}-\mbox{div} A(t,x,u,\nabla u)=\int_0^t(t-s)^{-\gamma}u(s)^{q}\,ds,} &\displaystyle {t>0,x\in {\mathbb{R}^d},}\\
\displaystyle u\geq0, &\displaystyle {t>0,x\in {\mathbb{R}^d},}\\
\displaystyle{u(0,x)= u_0(x),\qquad\qquad}&\displaystyle{x\in {\mathbb{R}^d},}
\end{array}
\right.
\]
where $q>1,  0<\gamma<1,  d \geq 1,  u_0\in L_{loc}^2(\mathbb{R}^d)$,  and $A:\mathbb{R}_+^*\times\mathbb{R}^d\times \mathbb{R}_+\times \mathbb{R}^d\longrightarrow \mathbb{R}^d$ is a Carath\'eodory function,  where one assumes the existence of $c_1,c_2>0$ and $p>1$ such that
\begin{equation}\label{condition1}
\left\{\begin{array}{l}
\displaystyle (A(t,x,u,w),w)\geq c_1\,|w|^p,\\
{}\\
\displaystyle |A(t,x,u,w)|\leq c_2\,|w|^{p-1},
\end{array}
\right.
\end{equation}
for all $(t,x,u,w)\in \mathbb{R}_+^*\times\mathbb{R}^d\times \mathbb{R}_+\times \mathbb{R}^d$.  For the study of the non-existence of global solution,  the following condition is needed
\begin{equation}\label{condition2}
 (A(t,x,u,w),w)\geq c_3\,|A(t,x,u,w)|^{\frac{p}{p-1}},\quad \mbox{for all} \,(t,x,u,w)\in \mathbb{R}_+^*\times\mathbb{R}^d\times \mathbb{R}_+\times \mathbb{R}^d,
\end{equation}
which is a general case of (\ref{condition1}).
\end{Remark}
Some words on the structure of the paper: Some definitions and properties on the
fractional integrals and derivatives are recalled in  Section \ref{sec2}.  Section \ref{sec3}, is reserved to the  proof of the main result (Theorem \ref{theo2}).

\section{Preliminaries}\label{sec2}

\begin{Definition}
A function $\mathcal{A} :[a,b]\rightarrow\mathbb{R}$, $-\infty<a<b<\infty$, is said to be absolutely continuous if and only if there exists $\psi\in L^1(a,b)$ such that
\[
\mathcal{A} (t)=\mathcal{A} (a)+\int_{a}^t \psi(s)\,ds.
\]
$AC[a,b]$ denotes the space of these functions.  Moreover, 
\[
AC^2[a,b]:=\left\{\varphi:[a,b]\rightarrow\mathbb{R}\;\hbox{/}\;\varphi'\in
AC[a,b]\right\}.
\]
\end{Definition}

\begin{Definition}
Let $f\in L^1(c,d)$, $-\infty<c<d<\infty$. The Riemann-Liouville left-and right-sided fractional integrals are defined by
\begin{equation}\label{I1}
I^\alpha_{c|t}f(t):=\frac{1}{\Gamma(\alpha)}\int_{c}^t(t-s)^{-(1-\alpha)}f(s)\,ds, \quad t>c, \; \alpha\in(0,1),
\end{equation}
and
\begin{equation}\label{I2}
I^\alpha_{t|d}f(t):=\frac{1}{\Gamma(\alpha)}\int_t^{d}(s-t)^{-(1-\alpha)}f(s)\,ds, \quad t<d, \; \alpha\in(0,1)
\end{equation}
where $\Gamma$ is the Euler gamma function.
\end{Definition}

\begin{Definition}
Let $f\in AC[c,d]$, $-\infty<c<d<\infty$. The Riemann-Liouville  left-and right-sided fractional derivatives  are defined by
\begin{equation}\label{}
D^\alpha_{c|t}f(t):=\frac{d}{dt}I^{1-\alpha}_{c|t}f(t)=\frac{1}{\Gamma(1-\alpha)}\frac{d}{dt}\int_{c}^t(t-s)^{-\alpha}f(s)\,ds, \quad t>c,  \; \alpha\in(0,1),
\end{equation}
and
\begin{equation}\label{}
D^\alpha_{t|d}f(t):=-\frac{d}{dt}I^{1-\alpha}_{t|d}f(t)=-\frac{1}{\Gamma(1-\alpha)}\frac{d}{dt}\int_t^{d}(s-t)^{-\alpha}f(s)\,ds, \quad t<d, \; \alpha\in(0,1).
\end{equation}
\end{Definition}

\begin{Proposition}\cite[(2.64), p. 46]{SKM}\\
Let $\alpha\in(0,1)$ and $-\infty<c<d<\infty$. The fractional integration by parts formula
\begin{equation}\label{IP}
\int_{c}^{d}f(t)D^\alpha_{c|t}g(t)\,dt \;=\; \int_{c}^{d}
g(t)D^\alpha_{t|d}f(t)\,dt,
\end{equation}
is satisfied for every $f\in I^\alpha_{t|d}(L^p(c,d))$, $g\in I^\alpha_{c|t}(L^q(c,d))$ such that $\frac{1}{p}+\frac{1}{q}\leq 1+\alpha$, $p,q>1$, where
\[
I^\alpha_{c|t}(L^q(0,T)):=\left\{f= I^\alpha_{c|t}h,\,\, h\in L^q(c,d)\right\},
\]
and
\[
I^\alpha_{t|d}(L^p(c,d)):=\left\{f= I^\alpha_{t|d}h,\,\, h\in L^p(c,d)\right\}.
\]
\end{Proposition}
\begin{Remark}
A simple sufficient condition for functions $f$ and $g$ to satisfy (\ref{IP}) is that $f,g\in C[c,d],$ such that
$D^\alpha_{t|d}f(t),D^\alpha_{c|t}g(t)$ exist at every point $t\in[c,d]$ and are continuous.
\end{Remark}

\begin{Proposition}\cite[Chapter~1]{SKM}\\
For $0<\alpha<1$, $-\infty<c<d<\infty$, we have the following identities
\begin{equation}\label{I3}
    D^\alpha_{c|t}I^\alpha_{c|t}\varphi(t)=\varphi(t),\,\hbox{a.e. $t\in(c,d)$}, \quad\hbox{for all}\,\varphi\in L^r(c,d), 1\leq r\leq\infty,
\end{equation}
and
\begin{equation}\label{I4}
   -\frac{d}{dt} \, D^\alpha_{t|d}\varphi=D^{1+\alpha}_{t|d}\varphi,\quad \varphi\in AC^2[c,d],
\end{equation}
where $D:=\frac{d}{dt}$.
\end{Proposition}

Given $T>0$, let us define the function $w_1$ by
\begin{equation}\label{w1}
\displaystyle w_1(t)=\left(1-t/T\right)^\sigma,\quad  t \in [0, T], \; \sigma\gg1
\end{equation}
The following properties concerning the functions $w_1$ will be used later on.
\begin{Lemma}\cite[(2.45), p. 40]{SKM}\label{L1}\\
Let $T>0$, $0<\alpha<1$. For all $t\in[0,T]$, we have
\begin{equation}\label{P1}
D_{t|T}^\alpha
w_1(t)=\frac{\Gamma(\sigma+1)}{
\Gamma(\sigma+1-\alpha)}T^{-\alpha}(1-t/T)^{\sigma-\alpha},
\end{equation}
and
\begin{equation}\label{P3}
D_{t|T}^{1+\alpha}
w_1(t)=\frac{\Gamma(\sigma+1)}{
\Gamma(\sigma-\alpha)}T^{-(1+\alpha)}(1-t/T)^{\sigma-\alpha-1}.
\end{equation}
\end{Lemma}
Any constant will be denoted by $C$.

\section{Global nonexistence}\label{sec3}
\setcounter{equation}{0}
In this section,  the proof of Theorem \ref{theo2} will be presented.  \\

We set 
$$Q_T=(0,T)\times \mathbb{R}^{d},\qquad\mathcal{B}_{T,B}=(0,T)\times \Omega(B),$$
where $\Omega(B):=\{x\in\mathbb{R}^n;\;|x|< 2B\}$, for $T,B>0$.

\begin{Definition}[Weak solution]\label{definitionweak}${}$\\
Let $0<\gamma<1$, and $T>0.$ A function
$$u\in C([0,T],L^{2}_{loc}(\mathbb{R}^d))\cap L^{2q}((0,T),L^{2q}_{loc}(\mathbb{R}^d))\cap L^{p}((0,T),W^{1,p}_{loc}(\mathbb{R}^d))\cap W^{1,2}((0,T),L^{2}_{loc}(\mathbb{R}^d)),$$
 is said to be a weak solution of (\ref{eq1}) if $u\geq0$, $u_{|_{t=0}}=u_0$ and the following formulation
\begin{eqnarray}\label{weaksolution}
\Gamma(\delta)\int_{Q_T}I^{\delta}_{0|t}(u^{q})\varphi\,dx\,dt= \int_{Q_T}|\nabla u|^{p-2}\nabla u\nabla \varphi\,dx\,dt+\int_{Q_T}u_t\varphi\,dx\,dt,
\end{eqnarray}
holds for all $\varphi\in  L^{2}(0,\infty;L^{2}(\mathbb{R}^d))\cap L^{p}(0;\infty,W^{1,p}(\mathbb{R}^d))$ such that supp$\varphi\subset[0,T]\times\mathbb{R}^d$ is compact, where $\delta=1-\gamma$. We denote the lifespan for the weak solution by
$$T_w(u_0):=\sup\{T\in(0,\infty];\,\,\hbox{there exists a unique weak solution u to (\ref{eq1})}\}.$$
Moreover, if $T>0$ can be arbitrary chosen, i.e. $T_w(u_0)=\infty$, then $u$ is called a global weak solution of (\ref{eq1}).
\end{Definition}

\begin{proof} Let $u \geq 0$ be a global weak solution of $(\ref{eq1})$, then
\[
\Gamma(\delta)\int_{Q_T}I^{\delta}_{0|t}(u^{q})\varphi\,dx\,dt= \int_{Q_T}|\nabla u|^{p-2}\nabla u\nabla \varphi\,dx\,dt+\int_{Q_T}u_t\varphi\,dx\,dt,
\]
for all $T>0$ and all $\varphi\in  L^{2}(0,\infty;L^{2}(\mathbb{R}^d))\cap L^{p}(0;\infty,W^{1,p}(\mathbb{R}^d))$ such that supp$\varphi\subset[0,T]\times\mathbb{R}^d$ is compact, where $\delta=1-\gamma$. Let
\begin{eqnarray}\label{testfunction}
\varphi(x,t)=D^\delta_{t|T}\left(\tilde{\varphi}(x,t)\right):=D^\delta_{t|T}\left(\varphi^\ell_1(x)\varphi_2(t)\right)\;\nonumber\\
\mbox{with}\quad\varphi_1(x):=\Phi\left({|x|}/{B}\right),\;\varphi_2(t):=\left(1-\frac{t}{T}\right)^\eta_+,
\end{eqnarray}
 where $\ell,\eta\gg1$ and $\Phi\in C^\infty(\mathbb{R}_+)$ be 
\[
\Phi(s)=\left\{\begin {array}{ll}\displaystyle{1},&\displaystyle{\quad 0\leq s\leq 1,}\\
\displaystyle{\searrow}&\displaystyle{\quad 1\leq s\leq 2,}\\
\displaystyle{0},&\displaystyle{\quad s\geq 2.}
\end {array}\right.
\]\\
Then,  using the integration-by-parts formula, we get
\begin{eqnarray*}
\Gamma(\delta)\int_{\mathcal{B}_{T,B}}I^{\delta}_{0|t}(u^{q})D^\delta_{t|T}\tilde{\varphi}\,dx\,dt&+&\int_{\Omega(B)}u_0(x)D^\delta_{t|T}\tilde{\varphi}(0,x)\,dx\\
&=&\int_{\mathcal{B}_{T,B}}|\nabla u|^{p-2}\nabla u\nabla \varphi\,dx\,dt-\int_{\mathcal{B}_{T,B}}u\;\partial_t\varphi\,dx\,dt.
\end{eqnarray*}
From $(\ref{IP})$ and Lemma \ref{L1}, we conclude that
\begin{eqnarray*}
\int_{\mathcal{B}_{T,B}}D^\delta_{0|t}I^{\delta}_{0|t}(u^{q})\tilde{\varphi}\,dx\,dt&+&C\;T^{-\delta}\int_{\Omega(B)}u_0(x)\varphi^\ell_1(x)\,dx\\
&=&C\int_{\mathcal{B}_{T,B}}|\nabla u|^{p-2}\nabla u\nabla \varphi\,dx\,dt-C\int_{\mathcal{B}_{T,B}}u\;\partial_t\varphi\,dx\,dt.
\end{eqnarray*}
Moreover, using (\ref{I3}), $u_0\geq0$ and $u\geq0$, it follows
\begin{equation}\label{estiA}
\int_{\mathcal{B}_{T,B}}u^{q}\tilde{\varphi}\,dx\,dt\leq C\int_{\mathcal{B}_{T,B}}|\nabla u|^{p-1}|\nabla \varphi|\,dx\,dt+C\int_{\mathcal{B}_{T,B}} u|\partial_t\varphi|\,dx\,dt=:I_1+I_2.
\end{equation}
 Next, by introducing $\tilde{\varphi}^{1/q}\tilde{\varphi}^{-1/q}$ in $I_2$ and applying the following Young's inequality
\[
AB\leq\frac{1}{8}A^q+C(q)B^{q'},\quad A\geq0,\;B\geq0,\;q+q'=q q',
\]
we get
\begin{equation}\label{conditionI1}
I_2\leq\frac{1}{8}\int_{\mathcal{B}_{T,B}}u^{q}\tilde{\varphi}\,dx\,dt+C\int_{\mathcal{B}_{T,B}}\varphi_1^\ell\varphi_2^{-\frac{1}{q-1}} \left(D^{1+\delta}_{t|T}\varphi_2\right)^{q'}\,dx\,dt.
\end{equation}
In order to obtain a similar estimation on $I_1$,  let $\alpha<0$ be an auxiliary constant such that $\alpha>\max\{-1,1-p,1-\frac{q}{p-1}\}$ and let
$$u_\epsilon(x,t)=u(x,t)+\epsilon,\quad\epsilon>0.$$
By taking $\varphi_\epsilon(x,t)=u_\epsilon^\alpha(x,t)\varphi(x,t)$ as a test function where $\varphi$ is given in (\ref{testfunction}), and using the fact that $u$ is a weak solution, we obtain
$$
\Gamma(\delta)\int_{\mathcal{B}_{T,B}}I^{\delta}_{0|t}(u^{q})u_\epsilon^\alpha\varphi\,dx\,dt= \int_{\mathcal{B}_{T,B}}|\nabla u|^{p-2}\nabla u\nabla(u_\epsilon^\alpha\varphi)\,dx\,dt+\int_{\mathcal{B}_{T,B}}u_t\,u_\epsilon^\alpha\varphi\,dx\,dt.
$$
Using the integration-by-parts formula,  we get
%$$$$
\begin{eqnarray*}
\Gamma(\delta)\int_{\mathcal{B}_{T,B}}I^{\delta}_{0|t}(u^{q})u_\epsilon^\alpha\varphi\,dx\,dt&=&\int_{\mathcal{B}_{T,B}}|\nabla u|^{p-2}\nabla u\nabla(u_\epsilon^\alpha\varphi)\,dx\,dt\\
&-&\int_{\mathcal{B}_{T,B}}u\partial_t(u_\epsilon^\alpha\varphi)\,dx\,dt -\displaystyle\int_{\Omega(B)}u_0(x)u_\epsilon^\alpha(x,0)\varphi(x,0)\,dx.
\end{eqnarray*}
Then,  as
$$
\nabla(u_\epsilon^\alpha\varphi)=\alpha u_\epsilon^{\alpha-1}\varphi\nabla u+u_\epsilon^\alpha\nabla\varphi\quad\mbox{and}\quad \partial_t(u_\epsilon^\alpha\varphi)=\alpha u_\epsilon^{\alpha-1}\varphi\partial_t u+u_\epsilon^\alpha\partial_t\varphi,
$$
it comes
\begin{eqnarray}\label{Eq1}
&{}&\Gamma(\delta)\int_{\mathcal{B}_{T,B}}I^{\delta}_{0|t}(u^{q})u_\epsilon^\alpha\varphi\,dx\,dt\nonumber\\
&{}&=\alpha \int_{\mathcal{B}_{T,B}}|\nabla u|^{p} u_\epsilon^{\alpha-1}\varphi\,dx\,dt + \int_{\mathcal{B}_{T,B}}|\nabla u|^{p-2}\nabla u\nabla\varphi\,u_\epsilon^\alpha\,dx\,dt\nonumber\\
&{}&\quad -\;\alpha\int_{\mathcal{B}_{T,B}}uu_{\epsilon}^{\alpha-1}\varphi\partial_tu\,dx\,dt-\int_{\mathcal{B}_{T,B}}uu_{\epsilon}^{\alpha}\partial_t\varphi\,dx\,dt
-\int_{\Omega(B)}u_0(x)u_\epsilon^\alpha(x,0)\varphi(x,0)\,dx.\quad\qquad
\end{eqnarray}
Whereupon,
\begin{eqnarray*}
J_1&:=&\int_{\mathcal{B}_{T,B}}uu_{\epsilon}^{\alpha-1}\varphi\partial_tu\,dx\,dt\\
&=&\int_{\mathcal{B}_{T,B}}(u_\epsilon-\epsilon)u_{\epsilon}^{\alpha-1}\varphi\partial_tu\,dx\,dt\\
&=&\int_{\mathcal{B}_{T,B}}u_{\epsilon}^{\alpha}\varphi\partial_tu_\epsilon\,dx\,dt-\epsilon\int_{\mathcal{B}_{T,B}}u_{\epsilon}^{\alpha-1}\varphi\partial_tu_\epsilon\,dx\,dt\\
&=&\frac{1}{\alpha+1}\int_{\mathcal{B}_{T,B}}\partial_t(u_{\epsilon}^{\alpha+1})\varphi\,dx\,dt-\frac{\epsilon}{\alpha}\int_{\mathcal{B}_{T,B}}\partial_t(u_\epsilon^{\alpha})\varphi\,dx\,dt\\
&=&-\frac{1}{\alpha+1}\int_{\mathcal{B}_{T,B}}u_{\epsilon}^{\alpha+1}\partial_t\varphi\,dx\,dt-\frac{1}{\alpha+1}\int_{\Omega(B)}u_{\epsilon}^{\alpha+1}(x,0)\varphi(x,0)\,dx\\
&{}&+\frac{\epsilon}{\alpha}\int_{\mathcal{B}_{T,B}}u_\epsilon^{\alpha}\partial_t\varphi\,dx\,dt+\frac{\epsilon}{\alpha}\int_{\Omega(B)}u_\epsilon^{\alpha}(x,0)\varphi(x,0)\,dx.
\end{eqnarray*}
Similarly,
\begin{eqnarray*}
J_2:&=&\int_{\mathcal{B}_{T,B}}uu_{\epsilon}^{\alpha}\partial_t\varphi\,dx\,dt\\
&=&\int_{\mathcal{B}_{T,B}}(u_\epsilon-\epsilon)u_{\epsilon}^{\alpha}\partial_t\varphi\,dx\,dt\\
&=&\int_{\mathcal{B}_{T,B}}u_{\epsilon}^{\alpha+1}\partial_t\varphi\,dx\,dt-\epsilon\int_{\mathcal{B}_{T,B}}u_{\epsilon}^{\alpha}\partial_t\varphi\,dx\,dt,
\end{eqnarray*}
and
\begin{eqnarray*}
J_3:&=&\int_{\Omega(B)}u_0(x)u_\epsilon^\alpha(x,0)\varphi(x,0)\,dx\\
&=&\int_{\Omega(B)}(u_\epsilon(x,0)-\epsilon)u_\epsilon^\alpha(x,0)\varphi(x,0)\,dx\\
&=&\int_{\Omega(B)}u_\epsilon^{\alpha+1}(x,0)\varphi(x,0)\,dx-\epsilon\int_{\Omega(B)}u_\epsilon^\alpha(x,0)\varphi(x,0)\,dx.
\end{eqnarray*}
By $J_1$,$J_2$, and $J_3$, it follows from (\ref{Eq1}) that
\begin{eqnarray*}
&&\Gamma(\delta)\int_{\mathcal{B}_{T,B}}I^{\delta}_{0|t}(u^{q})u_\epsilon^\alpha\varphi\,dx\,dt=\alpha \int_{\mathcal{B}_{T,B}}|\nabla u|^{p} u_\epsilon^{\alpha-1}\varphi\,dx\,dt +\\
&&\int_{\mathcal{B}_{T,B}}|\nabla u|^{p-2}\nabla u\nabla\varphi\,u_\epsilon^\alpha\,dx\,dt-\frac{1}{\alpha+1}\int_{\mathcal{B}_{T,B}}u_{\epsilon}^{\alpha+1}\partial_t\varphi\,dx\,dt\\
&&-\frac{1}{\alpha+1}\int_{\Omega(B)}u_{\epsilon}^{\alpha+1}(x,0)\varphi(x,0)\,dx,
\end{eqnarray*}
then
\begin{eqnarray*}
&{}& \Gamma(\delta)\int_{\mathcal{B}_{T,B}}I^{\delta}_{0|t}(u^{q})u_\epsilon^\alpha\varphi\,dx\,dt+|\alpha| \int_{\mathcal{B}_{T,B}}|\nabla u|^{p} u_\epsilon^{\alpha-1}\varphi\,dx\,dt \\
&{}&= \int_{\mathcal{B}_{T,B}}|\nabla u|^{p-2}\nabla u\nabla\varphi\,u_\epsilon^\alpha\,dx\,dt-\frac{1}{\alpha+1}\int_{\mathcal{B}_{T,B}}u_{\epsilon}^{\alpha+1}\partial_t\varphi\,dx\,dt\\
&&\quad-\frac{1}{\alpha+1}\int_{\Omega(B)}u_{\epsilon}^{\alpha+1}(x,0)\varphi(x,0)\,dx\\
&{}&\leq \int_{\mathcal{B}_{T,B}}|\nabla u|^{p-1}|\nabla\varphi |u_\epsilon^\alpha\,dx\,dt+\frac{1}{\alpha+1}\int_{\mathcal{B}_{T,B}}u_{\epsilon}^{\alpha+1}|\partial_t\varphi|\,dx\,dt\\
&{}&= \int_{\mathcal{B}_{T,B}}(|\nabla u|^{p-1}u_\epsilon^{\frac{(\alpha-1)(p-1)}{p}}\varphi^{\frac{p-1}{p}})(|\nabla\varphi |u_\epsilon^{\alpha-\frac{(\alpha-1)(p-1)}{p}}\varphi^{-\frac{p-1}{p}})\,dx\,dt\\
&&\quad+\frac{1}{\alpha+1}\int_{\mathcal{B}_{T,B}}u_{\epsilon}^{\alpha+1}|\partial_t\varphi|\,dx\,dt\\
&{}&\leq \frac{|\alpha|}{2} \int_{\mathcal{B}_{T,B}}|\nabla u|^{p} u_\epsilon^{\alpha-1}\varphi\,dx\,dt+C(\alpha) \int_{\mathcal{B}_{T,B}}u_{\epsilon}^{p-1+\alpha}|\nabla\varphi|^p\varphi^{1-p}\,dx\,dt \\ &&\quad+\frac{1}{\alpha+1}\int_{\mathcal{B}_{T,B}}u_{\epsilon}^{\alpha+1}|\partial_t\varphi|\,dx\,dt,
\end{eqnarray*}
where Young's inequality
\begin{equation}\label{young1}
AB\leq\frac{|\alpha|}{2}\;A^{\frac{p}{p-1}}+C(\alpha)\;B^p,\qquad A\geq0,\;B\geq0
\end{equation}
has been used.   
Consequently,  as $u \geq 0$,  we conclude that
$$
\int_{\mathcal{B}_{T,B}}|\nabla u|^{p} u_\epsilon^{\alpha-1}\varphi\,dx\,dt \leq C\int_{\mathcal{B}_{T,B}}u_{\epsilon}^{p-1+\alpha}|\nabla\varphi|^p\varphi^{1-p}\,dx\,dt +C \int_{\mathcal{B}_{T,B}}u_{\epsilon}^{\alpha+1}|\partial_t\varphi|\,dx\,dt.
$$
Young's inequality and  the last inequality,  allow to get
\begin{eqnarray*}
&{}&\int_{\mathcal{B}_{T,B}}|\nabla u|^{p-1} |\nabla\varphi|\,dx\,dt\\
&&= \int_{\mathcal{B}_{T,B}}\left(|\nabla u|^{p-1}u_\epsilon^{\frac{(\alpha-1)(p-1)}{p}}\varphi^{\frac{p-1}{p}}\right)\left(u_\epsilon^{\frac{(1-\alpha)(p-1)}{p}}\varphi^{\frac{1-p}{p}} |\nabla\varphi|\right)\,dx\,dt\\
&&\leq \int_{\mathcal{B}_{T,B}}|\nabla u|^{p} u_\epsilon^{\alpha-1}\varphi\,dx\,dt+ C\int_{\mathcal{B}_{T,B}}u_\epsilon^{(1-\alpha)(p-1)}\varphi^{1-p}|\nabla\varphi|^p\,dx\,dt\\
&&\leq C\int_{\mathcal{B}_{T,B}}u_{\epsilon}^{p-1+\alpha}|\nabla\varphi|^p\varphi^{1-p}\,dx\,dt +C \int_{\mathcal{B}_{T,B}}u_{\epsilon}^{\alpha+1}|\partial_t\varphi|\,dx\,dt\\
&{}&\quad+\;C \int_{\mathcal{B}_{T,B}}u_\epsilon^{(1-\alpha)(p-1)}\varphi^{1-p}|\nabla\varphi|^p\,dx\,dt.
\end{eqnarray*}
Applying Fatou's  and Lebesgue's theorems,  as $\epsilon\rightarrow0$,  we get \\
\begin{eqnarray*}
\displaystyle
I_1  &\leq& C\int_{\mathcal{B}_{T,B}}u^{p-1+\alpha}|\nabla\varphi|^p\varphi^{1-p}\,dx\,dt + C \int_{\mathcal{B}_{T,B}}u^{\alpha+1}|\partial_t\varphi|\,dx\,dt\\
& &+\,C \int_{\mathcal{B}_{T,B}}u^{(1-\alpha)(p-1)}\varphi^{1-p}|\nabla\varphi|^p\,dx\,dt=: K_1+K_2+K_3.
\end{eqnarray*}
Now, we use Young's inequality,   $\nabla(\varphi_1^\ell)=\ell\varphi_1^{\ell-1}\nabla\varphi_1$,  and the conditions that $q>\max\{1,p-1\}$, $\alpha<0$ to get estimations of $K_1$, $K_2$ and $K_3$;  we have
\begin{eqnarray*}
K_1&=&\int_{\mathcal{B}_{T,B}}(u^{p-1+\alpha}\tilde{\varphi}^{(p-1+\alpha)/q})(C\;\tilde{\varphi}^{-(p-1+\alpha)/q}|\nabla\varphi|^p\varphi^{1-p})\,dx\,dt\\
&\leq&\frac{1}{8}\int_{\mathcal{B}_{T,B}}u^{q}\tilde{\varphi}\,dx\,dt\\
&&+\,C\int_{\mathcal{B}_{T,B}}\varphi_1^{\ell-\frac{pq}{q-p+1-\alpha}}|\nabla\varphi_1|^{\frac{pq}{q-p+1-\alpha}}\varphi_2^{-\frac{p-1+\alpha}{q-p+1-\alpha}} \left(D^{\delta}_{t|T}\varphi_2\right)^{\frac{q}{q-p+1-\alpha}}\,dx\,dt,
\end{eqnarray*}
\begin{eqnarray*}
K_2&=&\int_{\mathcal{B}_{T,B}}(u^{\alpha+1}\tilde{\varphi}^{(\alpha+1)/q})(C\;\tilde{\varphi}^{-(\alpha+1)/q}|\partial_t\varphi|)\,dx\,dt \\
&\leq& \frac{1}{8}\int_{\mathcal{B}_{T,B}}u^{q}\tilde{\varphi}\,dx\,dt+C\int_{\mathcal{B}_{T,B}}\varphi_1^\ell\varphi_2^{-\frac{\alpha+1}{q-1-\alpha}} \left(D^{1+\delta}_{t|T}\varphi_2\right)^{\frac{q}{q-1-\alpha}}\,dx\,dt,
\end{eqnarray*}
and
\begin{eqnarray*}
K_3&=&\int_{\mathcal{B}_{T,B}}(u^{(1-\alpha)(p-1)}\tilde{\varphi}^{(1-\alpha)(p-1)/q})(C\;\tilde{\varphi}^{-(1-\alpha)(p-1)/q}\varphi^{1-p}|\nabla\varphi|^p)\,dx\,dt\\\
&\leq&\frac{1}{8}\int_{\mathcal{B}_{T,B}}u^{q}\tilde{\varphi}\,dx\,dt\\
&&+\,C\int_{\mathcal{B}_{T,B}}\varphi_1^{\ell-\frac{pq}{q-(1-\alpha)(p-1)}}|\nabla\varphi_1|^{\frac{pq}{q-(1-\alpha)(p-1)}}\varphi_2^{-\frac{(1-\alpha)(p-1)}{q-(1-\alpha)(p-1)}} \left(D^{\delta}_{t|T}\varphi_2\right)^{\frac{q}{q-(1-\alpha)(p-1)}}\,dx\,dt.
\end{eqnarray*}
Therefore, we conclude that
\begin{eqnarray*}
I_1&\leq& \frac{3}{8}\int_{\mathcal{B}_{T,B}}u^{q}\tilde{\varphi}\,dx\,dt\\
&+&C\int_{\mathcal{B}_{T,B}}\varphi_1^{\ell-\frac{pq}{q-p+1-\alpha}}|\nabla\varphi_1|^{\frac{pq}{q-p+1-\alpha}}\varphi_2^{-\frac{p-1+\alpha}{q-p+1-\alpha}} \left(D^{\delta}_{t|T}\varphi_2\right)^{\frac{q}{q-p+1-\alpha}}\,dx\,dt\\
&{}&+\;C\int_{\mathcal{B}_{T,B}}\varphi_1^\ell\varphi_2^{-\frac{\alpha+1}{q-1-\alpha}} \left(D^{1+\delta}_{t|T}\varphi_2\right)^{\frac{q}{q-1-\alpha}}\,dx\,dt\\
&{}&+\; C\int_{\mathcal{B}_{T,B}}\varphi_1^{\ell-\frac{pq}{q-(1-\alpha)(p-1)}}|\nabla\varphi_1|^{\frac{pq}{q-(1-\alpha)(p-1)}}\varphi_2^{-\frac{(1-\alpha)(p-1)}{q-(1-\alpha)(p-1)}} \left(D^{\delta}_{t|T}\varphi_2\right)^{\frac{q}{q-(1-\alpha)(p-1)}}\,dx\,dt.
\end{eqnarray*}
Using the estimates of $I_1$ and $I_2$ into (\ref{estiA}), we obtain
\begin{eqnarray}\label{2000}
&&\int_{\mathcal{B}_{T,B}}u^{q}\tilde{\varphi}\,dx\,dt\nonumber\\
&&\leq C\int_{\mathcal{B}_{T,B}}\varphi_1^\ell\varphi_2^{-\frac{1}{q-1}} \left(D^{1+\delta}_{t|T}\varphi_2\right)^{\frac{q}{q-1}}\,dx\,dt\nonumber\\
&&\quad +\,C\int_{\mathcal{B}_{T,B}}\varphi_1^{\ell-\frac{pq}{q-p+1-\alpha}}|\nabla\varphi_1|^{\frac{pq}{q-p+1-\alpha}}\varphi_2^{-\frac{p-1+\alpha}{q-p+1-\alpha}} \left(D^{\delta}_{t|T}\varphi_2\right)^{\frac{q}{q-p+1-\alpha}}\,dx\,dt\nonumber\\
&&\quad+\,C\int_{\mathcal{B}_{T,B}}\varphi_1^\ell\varphi_2^{-\frac{\alpha+1}{q-1-\alpha}} \left(D^{1+\delta}_{t|T}\varphi_2\right)^{\frac{q}{q-1-\alpha}}\,dx\,dt\nonumber\\
&&\quad+\,C\int_{\mathcal{B}_{T,B}}\varphi_1^{\ell-\frac{pq}{q-(1-\alpha)(p-1)}}|\nabla\varphi_1|^{\frac{pq}{q-(1-\alpha)(p-1)}}\varphi_2^{-\frac{(1-\alpha)(p-1)}{q-(1-\alpha)(p-1)}} \left(D^{\delta}_{t|T}\varphi_2\right)^{\frac{q}{q-(1-\alpha)(p-1)}}\,dx\,dt.\qquad
\end{eqnarray}
At this stage, we choose $B=T^{\theta}$, $\theta>0$. Taking $\alpha$ small enough and passing to  $s=T^{-1}t,$ $y=T^{-\theta}x$, we get from (\ref{2000}) that
\begin{equation}\label{newweaksolution5}
\int_0^T\int_{\Omega(T^{\theta})}|u|^{q}\tilde{\varphi}\,dx\,dt\leq C\,T^{-\frac{(\delta+1)q}{q-1}+d\theta+1}+ C\,T^{-\frac{(\delta+\theta p) q}{p-p+1}+d\theta+1}.
\end{equation}
If all exponents of $T$ are negative, by taking $T\longrightarrow+\infty$ and using the dominated convergence theorem, we conclude that $u=0$. In order to ensure the negativity of the exponents of $T$,  it is sufficient to require
$$q<\frac{1+d\theta}{(d\theta-\delta)_+}=:q_1(\theta)\quad\hbox{and}\quad q<\frac{(p-1)(1+d\theta)}{(d\theta-\theta p+1-\delta)_+}=:q_2(\theta).$$
which is equivalent to
\begin{equation}\label{q1q2}
q<\max_{\theta>0}\min\{q_1(\theta),q_2(\theta)\}.
\end{equation}
To take into consideration $q_1(\theta)$ and $q_2(\theta)$, we first look at $(d\theta-\delta)_+$ and $(d\theta-\theta p+1-\delta)_+$ and try to compare them in terms of $\theta$, i.e,  to compare between $\delta/d$ and $(1-\delta)/((p-d)_+)$;  this requires to study two cases: $p\geq d/\delta$ and $p<d/\delta$.\\

\noindent{\bf\underline{$a)$ If $p\geq d/\delta$:}}  In this case we have $(1-\delta)/(p-d)\leq \delta/d$. As $q_1(\theta)$ is non-increasing and $q_2(\theta)$ is non-decreasing as a function of $\theta$, (\ref{q1q2}) can be read as
$$q<\max\left\{\max_{\theta\in[\frac{\delta}{d},+\infty[}q_1(\theta),\max_{\theta\in]0,\frac{1-\delta}{p-d}]}q_2(\theta)\right\}=+\infty.$$

\noindent{\bf\underline{$b)$ If $p< d/\delta$:}} Technical calculations lead us to to distinguish $3$ cases.\\

{\bf\underline{$i)$ Case $p> d$.}}  In this case we have $\delta/d<(1-\delta)/(p-d)$. So, (\ref{q1q2}) can be read as
$$q<\max\left\{\max_{\theta\in[\theta_0,+\infty[}q_1(\theta),\max_{\theta\in]0,\theta_0]}q_2(\theta)\right\}=q_1(\theta_0)=q_2(\theta_0)=q_c=\max\left\{\frac{p-1}{1-\delta},q_c\right\},$$
with
$$\theta_0:=\frac{\delta p+1-2\delta}{dp+p-2d}>0,$$
using that $q_1(\theta)$ is non-increasing and $q_2(\theta)$ is non-decreasing.\\

{\bf\underline{$ii)$ Case $2d/(d+1)<p\leq d$.}} In this case, (\ref{q1q2}) can be read as
\begin{equation}\label{d_0}
q<\max\left\{\max_{\theta\in[\theta_0,+\infty[}q_1(\theta),\max_{\theta\in]0,\theta_0]}q_2(\theta)\right\}.
\end{equation}
As $p_2(\theta)$ is non-decreasing when $p\geq d\delta$, and non-increasing when $p\leq d\delta$, we can see that (\ref{d_0}) is equivalent to
$$q<\max\left\{\frac{p-1}{1-\delta},q_c\right\}.$$

{\bf\underline{$iii)$ Case of $1<p\leq 2d/(d+1)$:}} In this case, (\ref{q1q2}) can be read as
$$
q<\max\left\{\max_{\theta\in]0,+\infty[}q_2(\theta)\right\}=q_2(0)=\frac{p-1}{1-\delta}=\max\left\{\frac{p-1}{1-\delta},q_c\right\}.
$$
Finally, to get similar results in the critical case
$$q=\max\left\{\frac{p-1}{1-\delta},q_c\right\},$$
we choose $B=R^{-\theta}T^{\theta},$ where $1\ll R<T$ is such that $T$ and $R$ do not go simultaneously to infinity. Moreover, due to the calculation made above,  a positive constant $D$ independent of $T$ exists such that
$$\int_0^\infty\int_{\mathbb{R}^d}|u|^q\,dx\,dt\leq D,$$
leading to
\begin{equation}\label{conditioninfini}
\int_0^T\int_{\Delta(B)}|u|^q\tilde{\varphi}\,dx\,dt\rightarrow0\quad\mbox{as}\;\;T\rightarrow\infty,
\end{equation}
where $\Delta(B):=\{x\in\mathbb{R}^d;\;B<|x|< 2B\}$. Repeating a similar calculation as in the subcritical case, $q<\max\{\frac{p-1}{1-\delta},q_c\}$, and using H\"older's inequality instead of Young's one in $K_1$ and $K_3$,  we get
$$
I_2\leq\frac{1}{3}\int_{\mathcal{B}_{T,B}}u^{q}\tilde{\varphi}\,dx\,dt+C\int_{\mathcal{B}_{T,B}}\varphi_1^\ell\varphi_2^{-1/(q-1)} \left(D^{1+\delta}_{t|T}\varphi_2\right)^{q'}\,dx\,dt,
$$
\begin{eqnarray*}
&{}&K_1 \leq C\left(\int_0^T\int_{\Delta(B)}u^{q}\tilde{\varphi}\,dx\,dt\right)^{\frac{p-1+\alpha}{q}}\\
&&\quad\qquad\left(\int_{\mathcal{B}_{T,B}}\varphi_1^{\frac{\ell(q-p+1-\alpha)-pq}{q-p+1-\alpha}}|\nabla\varphi_1|^{\frac{pq}{q-p+1-\alpha}}\varphi_2^{-\frac{p-1+\alpha}{q-p+1-\alpha}} \left(D^{\delta}_{t|T}\varphi_2\right)^{\frac{q}{q-p+1-\alpha}}\,dx\,dt\right)^{\frac{q-p+1-\alpha}{q}},
\end{eqnarray*}

\begin{eqnarray*}
K_2&\leq& \frac{1}{3}\int_{\mathcal{B}_{T,B}}u^{q}\tilde{\varphi}\,dx\,dt+C\int_{\mathcal{B}_{T,B}}\varphi_1^\ell\varphi_2^{-\frac{\alpha+1}{q-1-\alpha}} \left(D^{1+\delta}_{t|T}\varphi_2\right)^{\frac{q}{q-1-\alpha}}\,dx\,dt,
\end{eqnarray*}
and

\begin{eqnarray*}
&&K_3\leq C \left(\int_0^T\int_{\Delta(B)}u^{q}\tilde{\varphi}\,dx\,dt\right)^{\frac{(1-\alpha)(p-1)}{q}}\\
&&\left(\int_{\mathcal{B}_{T,B}}\varphi_1^{\frac{\ell(q-(1-\alpha)(p-1))-pq}{q-(1-\alpha)(p-1)}}|\nabla\varphi_1|^{\frac{pq}{q-p+1-\alpha}}\varphi_2^{-\frac{(1-\alpha)(p-1)}{q-(1-\alpha)(p-1)}} \left(D^{\delta}_{t|T}\varphi_2\right)^{\frac{q}{q-(1-\alpha)(p-1)}}\,dx\,dt\right)^{\frac{q-(1-\alpha)(p-1)}{q}}.
\end{eqnarray*}
Whereupon
\begin{eqnarray}\label{2001}
&&\int_{\mathcal{B}_{T,B}}u^{q}\tilde{\varphi}\,dx\,dt\nonumber\\
&&\leq\;C\int_{\mathcal{B}_{T,B}}\varphi_1^\ell\varphi_2^{-1/(q-1)} \left(D^{1+\delta}_{t|T}\varphi_2\right)^{\frac{q}{q-1}}\,dx\,dt\nonumber\\
&&+C \left(\int_0^T\int_{\Delta(B)}u^{q}\tilde{\varphi}\,dx\,dt\right)^{\frac{p-1}{q}}\nonumber\\
&&\qquad\left(\int_{\mathcal{B}_{T,B}}\varphi_1^{\frac{\ell(q-p+1)-pq}{q-p+1}}|\nabla\varphi_1|^{\frac{pq}{q-p+1}}\varphi_2^{-\frac{p-1}{q-p+1}} \left(D^{\delta}_{t|T}\varphi_2\right)^{\frac{q}{q-p+1}}\,dx\,dt\right)^{\frac{q-p+1}{q}}\nonumber\\
&&+C\int_{\mathcal{B}_{T,B}}\varphi_1^\ell\varphi_2^{-\frac{1}{q-1}} \left(D^{1+\delta}_{t|T}\varphi_2\right)^{\frac{q}{q-1}}\,dx\,dt\nonumber\\
&&+C \left(\int_0^T\int_{\Delta(B)}u^{q}\tilde{\varphi}\,dx\,dt\right)^{\frac{p-1}{q}}\nonumber\\
&&\qquad\left(\int_{\mathcal{B}_{T,B}}\varphi_1^{\frac{\ell(q-p+1)-pq}{q-p+1}}|\nabla\varphi_1|^{\frac{pq}{q-p+1}}\varphi_2^{-\frac{p-1}{q-p+1}} \left(D^{\delta}_{t|T}\varphi_2\right)^{\frac{q}{q-p+1}}\,dx\,dt\right)^{\frac{q-p+1}{q}}.
\end{eqnarray}
Taking into account that $q=\max\{\frac{p-1}{1-\delta},q_c\}$ and $s=T^{-1}t, \;y=R^{\theta}T^{-\theta}x$, we get
 $$\int_{\mathcal{B}_{T,B}}|u|^q\tilde{\varphi}\,dx\,dt \leq C\,R^{-d\theta}+C R^{-d\theta+\frac{pq\theta}{q-p+1}} \left(\int_0^T\int_{\Delta(B)}u^{q}\tilde{\varphi}\,dx\,dt\right)^{\frac{p-1}{q}}$$
Taking the limit when $T\rightarrow\infty$, and using $(\ref{conditioninfini})$,  we obtain
  $$\int_0^\infty\int_{\mathbb{R}^d}|u|^q\,dx\,dt \leq  C\,R^{-d\theta}.$$
 Finally, letting $R\rightarrow\infty$,  it comes  that $u=0$. The proof of Theorem $\ref{theo2}$ is complete.
\end{proof}

\end{document}